\documentclass[oneside,draft,12pt]{amsart}
\usepackage{amssymb, amscd}
\usepackage[all]{xy}

\theoremstyle{plain}
\newtheorem{theorem}{Theorem}[section]
\newtheorem{lemma}[theorem]{Lemma}
\newtheorem*{thm}{Theorem}
\newtheorem{proposition}[theorem]{Proposition}
\newtheorem*{prop}{Proposition}
\newtheorem{corollary}[theorem]{Corollary}

\theoremstyle{definition}

\newtheorem{definition}[theorem]{Definition}

\newtheorem{example}[theorem]{Example}

\theoremstyle{remark}
\newtheorem{remark}[theorem]{Remark}
\newtheorem{notation}[theorem]{Notation}
\newtheorem{question}[theorem]{Question}

\DeclareMathOperator{\Witt}{\textit{W}}

\DeclareMathOperator{\Ker}{Ker}

\DeclareMathOperator{\Hom}{Hom}

\DeclareMathOperator{\Gal}{Gal}

\DeclareMathOperator{\Image}{Im}
\DeclareMathOperator{\level}{s}

\DeclareMathOperator{\cd}{cd}

\newcommand{\mcG}{\ensuremath{\mathcal{G}}}

\newcommand{\bbF}{\ensuremath{\mathbb{F}}}
\newcommand{\bbN}{\ensuremath{\mathbb{N}}}
\newcommand{\bbR}{\ensuremath{\mathbb{R}}}
\newcommand{\bbQ}{\ensuremath{\mathbb{Q}}}
\newcommand{\bbC}{\ensuremath{\mathbb{C}}}
\newcommand{\bbZ}{\ensuremath{\mathbb{Z}}}
\newcommand{\bbS}{\ensuremath{\mathbb{S}}}
\newcommand{\Fbb}{\ensuremath{\mathbb{F}}}

\newcommand{\ra}{\ensuremath{\rightarrow}}

\begin{document}

\title{Field Theory and the Cohomology of Some Galois Groups}

\author{Alejandro Adem}
\address{Mathematics Department\\
         University of Wisconsin\\
         Madison, Wisconsin, 53706} 
\email{adem@math.wisc.edu}
\thanks{The first author was partially supported by an NSF grant, the
MPIM-Bonn, and the CRM-Barcelona.}

\author{Wenfeng Gao}
\address{4426 S. ${112}^{\text{th}}$ E.\ Ave. \\
         Tulsa, OK}
\email{wgao@sprintmail.com}

\author{Dikran B. Karagueuzian}
\address{Mathematics Department\\
         University of Wisconsin\\
         Madison, Wisconsin, 53706} 
\email{dikran@math.wisc.edu}
\thanks{The third author was 
partially supported by an NSF postdoctoral fellowship, CRM-Barcelona,
and MPIM-Bonn.}

\author{J\'an Min\'a\v{c}}
\address{Mathematics Department \\
         University of Western Ontario \\
         London, Ontario, Canada N6A 5B7}
\email{minac@math.uwo.ca.edu}
\thanks{The fourth author was partially supported by the NSERC and the
special Dean of Science fund at UWO. }


\begin{abstract}

We prove that two arithmetically significant extensions of a field $F$
coincide if and only if the Witt ring $\Witt F$ is a group ring
$\bbZ/n[G]$.  Furthermore, working modulo squares with Galois groups
which are $2$-groups, we establish a theorem analogous to
Hilbert's theorem 90 and show that an identity linking the
cohomological dimension of the Galois group of the quadratic closure of
$F$, the length of a filtration on a certain module over a
Galois group, and the
dimension over $\bbF_2$ of the square class group of the field holds for
a number of interesting families of fields.  Finally, we discuss the
cohomology of a particular Galois group in a topological
context. 

\end{abstract}

\maketitle

\begin{section}{Introduction}
\label{s:intro}

In recent years substantial breakthroughs have been made in the study of
the cohomology of absolute Galois groups, culminating in Voevodsky's
proof of the Milnor conjecture \cite{Voe}.  These results severely
restrict possible absolute Galois groups of fields and it seems
virtually certain that they will have purely field-theoretic
consequences.  Such results are not easily derived, however, and in
fact only a few theorems of this type have appeared (for some
examples the reader may consult \cite{GMi, MS1, Vill}). The goal of
this paper is to obtain some results in field theory as consequences
of Merkurjev's theorem \cite{Mer}.  The techniques we use are somewhat
varied; for example, we 
study the square class group of a field as a module for a Galois group
and relate its socle series to the $E^{1,1}_\infty$-term
of a spectral sequence. We are then able to obtain information on this term 
of the spectral sequence using techniques from the theory of binary
quadratic forms. In the study of particular examples our techniques
become more specialized, for example, at one point we make essential use
of the main theorem of local class field theory.

Throughout the paper, $F$ will denote a field whose characteristic is
not $2$. We write
$F^{(2)}$ for the field obtained by
adjoining all the square roots of elements of $F$. 
Now let $F^{\{3\}}= (F^{(2)})^{(2)}$
and $\Gal(F^{(2)}/F)=G_F^{[2]}$.
We introduce  
the following definition:
the \textit{V}-\textsl{group} of $F$ is the Galois group
$G_F^{\{3\}} = \Gal(F^{\{3\}}/F)$. 
In \cite{MS1} Min\'a\v{c} and Spira defined an
extension of $F$ whose Galois group is closely related to the Witt
ring of $F$. 
Let $F^{(3)}$ be the extension of $F^{(2)}$ obtained by adjoining the
square roots of elements $\alpha \in F^{(2)}$ such that
$F^{(2)}(\sqrt{\alpha})$ is Galois over $F$.  $F^{(3)}$ is a Galois
extension of $F$ and we recall
that the \textit{W}-\textsl{group} of $F$ is the Galois group
${\mathcal G}_F = \Gal(F^{(3)}/F)$.
Of course $F^{(3)} \subset F^{\{3\}}$, so that ${\mathcal G}_F$ is a
quotient of $G_F^{\{3\}}$.
Furthermore, it follows immediately from the
definitions that ${\mathcal G}_F = G_F^{\{3\}}$ if and only if
$F^{\{3\}} = F^{(3)}$.

In this paper we use tools from the cohomology of groups
to study these Galois groups and the field-theoretic information
they may contain.
Our first result is a characterization of \textit{C}-fields
(i.e. fields $F$ such that the Witt ring $\Witt F$ is isomorphic
to a group ring of the form $\bbZ/n\bbZ[G]$):

\begin{thm}{\rm \textbf{(\ref{t:C-field-implies})}}
 The \textit{W}-group of $F$ is the same as the \textit{V}-group of $F$
if and only if $F$ is a \textit{C}-field.
\end{thm}

The proof makes use of a theorem due to Merkurjev as well
as a description of the $E^{1,1}_{\infty}$ term of the mod 2
Lyndon-Hochschild-Serre spectral sequence associated to the
Frattini extension for ${\mathcal G}_F$.
We introduce a length invariant associated to $F$ as follows.
Let $J := \dot{F^{(2)}}/(\dot{F^{(2)}})^2$, regarded as a
$G_F^{[2]}$-module.  The socle series of $J$,
\[
0 = J_0 \subset J_1 \subset J_2 \subset \dots \subset J
\]
has finite length $l$ (i.e. $J_l = J$) if $\dot{F}/\dot{F}^2$ is finite.  
This length is an invariant of the field, which we will write as
$l(F)$.  

One of the key results used in our proof is 
\begin{prop}{\rm \textbf{(\ref{p:Mer})}}
In the mod $2$
Lyndon--Hochschild--Serre spectral sequence for the group
extension 
$1\to \Phi \rightarrow {\mathcal G}_F \rightarrow
G_F^{[2]}\to 1$, 
we have
$E_{\infty}^{1,1} \approx J_2/J_1$.
\end{prop}

Next we obtain a result reminiscent of Hilbert's Theorem 90,
describing the elements in the kernel of the differential $d_2$.
While the proof of this theorem depends on Proposition~\ref{p:Mer} and
hence indirectly on Merkurjev's theorem, we provide a combinatorial,
group-theoretic proof of ``half'' of the theorem in an appendix
(\S\ref{s:appendix}).

Let $\{[a_i] \mid i\in I\}$ be a basis for $\dot F/\dot F^2$,
which we identify with $H^1(G_F^{[2]},\bbF_2)$. Let $\{\sigma_i
\mid i\in I\}$ be a minimal set of generators of $G_F^{[2]}$ such that
${\sigma_i(\sqrt{a_j})}/{\sqrt{a_j}}= (-1)^{\delta_{i,j}}$, where
$i,j\in I$ and $\delta_{i,j}$ is the Kronecker delta function.
One may view the $a_i$ as a dual basis to the $\sigma_i$. 
\begin{thm}{\rm \textbf{(\ref{p:J90})}}
If $\left\{[j_{i_k}] \mid 1\leq k\leq l\right\}$ is a  finite
set of elements in $J_1$, then 
\[
    \sum_{k=0}^l[a_{i_k}]d_2^{0,1}\left([j_{i_k}]\right) =0
\]
if and only if there exists an element $[j]\in J_2$ such that
\begin{enumerate}
\item $[{\sigma_{i_k}(j)}/{j}] =[j_{i_k}]$ for each $k\in\{1,\ldots,l\}$ and 
\item $[\sigma_u(j)] =[j] $ for each $u \notin \{i_1,\ldots,i_l\}$. 
\end{enumerate}
\end{thm}

One should observe that this result is a practical criterion
for the existence of elements $[j]$ from $J_2$ with prescribed images
under the action of the Galois group $G_F^{[2]}$. We illustrate this
in a detailed way in a special example when the $W$--group is
the universal $W$--group on two generators (see example 4.2).
Moreover this theorem fully describes the $G_F^{[2]}$--module
$J_2$.

Let $G_q := \Gal(F_q/F)$, where $F_q$ is a quadratic closure of $F$.
Our next result provides a formula for the cohomological
dimension of $G_q$ for certain non-formally real fields in terms
of the length function $l(F)$:

\begin{thm}{\rm \textbf{(\ref{p:F=C-field}, \ref{e:GF=Wn},
\ref{t:F=local})}}
Let $F$ be a field of which is not formally real, and assume that 
$|\dot F/\dot F^2|=2^n$. Then if $F$ is a \textit{C}-field, a local
field, or a field such that $G_q$ is free, we have
$l(F)+ \cd(G_q)=n+1$.
\end{thm}

Finally in our last section we make an explicit analysis
of the mod 2 cohomology of the \textit{V}-group arising when $G_q$
is the free pro-2-group on 2 generators (denoted $V(2)$, a group
of order $2^7$). We obtain its
Poincar\'e series as well as some rather interesting topological
information.

\begin{prop}{\rm \textbf{(\ref{p:V2})}}
The Poincar\'e series for $H^*(V(2),\bbF_2)$ is 
\[
\frac{1-x+x^2}{(1-x)^3(1-x^2)^2},
\]
which expands to $1 + 2 x + 6 x^2  + 11 x^3  + 22 x^4  + 36 x^5 + 60
x^6 + 90 x^7 + 135 x^8 + 190 x^9 + 266x^{10} +
\cdots$.  
\end{prop}

Our methods combine techniques from field theory with methods
from the cohomology of groups. 
For background information on Galois cohomology, group cohomology,
quadratic form theory, and local class field theory, we refer the reader
to \cite{S-GC}, \cite{AM,Benson}, \cite{Lam1, Lam2}, and
\cite{AT,CF} respectively.  We have endeavored to provide
specific references throughout the paper as often as possible.

This paper is organized as
follows: in \S2 we describe
preliminary definitions and basic facts about the Galois groups
considered here; in \S3 we prove the characterization of
\textit{C}-fields and introduce the length of $F$; in \S4 we prove
our analogue of Hilbert's Theorem 90; in \S5 we discuss our results
relating $l(F)$ with $\cd(G_q)$ for examples such as local fields;
and finally in \S6 we provide a fairly complete cohomological analysis
of the group $V(2)$ and of a $5$-dimensional crystallographic group
associated to it.

A few comments about notation may be helpful.
We shall have several uses for the cyclic group of order two, and thus
as many notations for it.  Generally, we denote this group by
$\bbZ/2\bbZ$, but when we regard it as the group of square roots of
unity in a field, we write it $\mu_2$, and when we consider the field of
two elements, used as coefficients for cohomology, we write $\bbF_2$.
Generally, we do not specify the coefficients for our cohomology
theories, and in such cases $\bbF_2$-coefficients are to be understood.
When we use other coefficient
rings, they will be specified.

\end{section}

\begin{section}{Preliminaries}
\label{s:preliminaries}

Let $F$ denote a field of characteristic different from two, and
denote by $F^{(2)}$ the field obtained by
adjoining all the square roots of elements of $F$. 
Now let $F^{\{3\}}= (F^{(2)})^{(2)}$.
 From 
Kummer theory it is easy to see that
$F^{(2)}$ is a Galois extension of $F$, and if $B$ is a basis of $\dot
F / \dot F^2$, then the Galois group of $F^{(2)}$ over $F$ is 
$\prod_B {\mathbb Z}/2\bbZ$. 
We denote this Galois group by $G_F^{[2]}$.
A slightly less obvious fact is 
\begin{lemma}
\label{l:F3/F-Galois}
$F^{\{3\}}$ is a Galois extension of $F$. 
\end{lemma}

\begin{proof}
Let $\sigma \colon F^{\{3\}} \rightarrow \bar{F}$ be an embedding of
$F^{\{3\}}$ into an algebraic closure $\bar{F}$ of $F$ containing
$F^{\{3\}}$.  It is enough to show that for each $\gamma \in
\dot{F^{(2)}}$, $\sigma(\sqrt{\gamma}) \in F^{\{3\}}$.  However,
$\sigma(\sqrt{\gamma})^2 = \sigma(\gamma)$ and therefore
$\sigma(\sqrt{\gamma}) = \pm \sigma(\sqrt{\gamma})$.  But
$\sigma(\gamma) \in \dot{F^{(2)}}$, since $F^{(2)}/F$ is Galois, so
the lemma follows.
\end{proof}

Using this we introduce  

\begin{definition}
The \textit{V}-\textsl{group} of $F$ is the Galois group
$G_F^{\{3\}} = \Gal(F^{\{3\}}/F)$. 
\end{definition}

In \cite{MS1} Min\'a\v{c} and Spira defined  the ``Witt closure'' of $F$, an
extension of $F$ whose Galois group is closely related to the Witt
ring of $F$. 
Let $F^{(3)}$ be the extension of $F^{(2)}$ obtained by adjoining the
square roots of elements $\alpha \in F^{(2)}$ such that
$F^{(2)}(\sqrt{\alpha})$ is Galois over $F$.  $F^{(3)}$ is a Galois
extension of $F$ and we recall

\begin{definition}
The \textit{W}-\textsl{group} of $F$ is the Galois group
${\mathcal G}_F = \Gal(F^{(3)}/F)$.
\end{definition}

Of course $F^{(3)} \subset F^{\{3\}}$, so that ${\mathcal G}_F$ is a
quotient of $G_F^{\{3\}}$. We will frequently make use of this fact
without explicit mention. Furthermore, it follows immediately from the
definitions that ${\mathcal G}_F = G_F^{\{3\}}$ if and only if
$F^{\{3\}} = F^{(3)}$. 

Now we introduce the Galois group of the quadratic
closure and one of its subgroups:

\begin{definition}
Let $G_q := \Gal(F_q/F)$, where $F_q$ is a quadratic closure of $F$,
and set $G^{(2)} := \Gal(F_q/F^{(2)})$. 
\end{definition}

Note that $G_q/G^{(2)}\cong G_F^{[2]}$ and 
that the \textit{V}-group is determined by $G_q$, namely if
$\Phi (G^{(2)})$ denotes the Frattini subgroup of $G^{(2)}$, then
from our definition it follows that
\[
G_F^{\{3\}} = G_q / \Phi (G^{(2)}).
\]  

 From the above we obtain that there is an extension of elementary 
abelian groups 
\[
1\rightarrow G_{F^{(2)}}^{[2]} \rightarrow G_F^{\{3\}} \rightarrow
G_{F}^{[2]}\rightarrow 1,
\]
where the action on the kernel can be highly non-trivial.
We will often abbreviate $G_{F^{(2)}}^{[2]}$ by simply using the symbol $A$.
 From the point of view of group theory, the \textit{V}-group of the field $F$
is simply the extension of $G_F^{[2]}$ obtained by taking the quotient
of $G_q$ by the
Frattini subgroup of $G^{(2)}$.

We will be especially interested in the case when $G_q$ is a free pro-2-group
on $n$ generators. From the analysis in \cite{MS1} we know that in this
case the \textit{W}-group maps onto any \textit{W}-group arising from a field
with $n$-dimensional group of square classes, hence it is referred to
as the universal \textit{W}-group $W(n)$. Similarly the \textit{V}-group arising from
this situation will enjoy an analogous universal property and we denote it
by $V(n)$. 

In this case the group $G^{(2)}$
can be identified with the Frattini subgroup of the free
pro-2-group $G_q$, hence it is a free pro-2-group of rank $2^n(n-1)+1$,
where $G_F^{[2]}\cong (\bbZ/2\bbZ)^n$. More explicitly we have
$G^{(2)}=G_q^2[G_q,G_q]$ (in fact, $[G_q,G_q] \subset G_q^2$, so
$G^{(2)} = G_q^2$).  Our next step is to ``abelianize'' this extension,
namely we factor out the Frattini subgroup of $G^{(2)}$ 
This yields the quotient $V(n)$ as an extension
\[
1\to A\to V(n)\to G_F^{[2]}\to 1. 
\]
Moreover we can identify the $G_F^{[2]}$-module $A$ explicitly using the methods
in \cite{Ku}; indeed $A\cong \Omega^2(\bbF_2)$, the second dimension
shift of the trivial module (for the definition of the ``Heller
translate'' $\Omega^2$ see \cite[p.~8, v.~I]{Benson}). Alternatively, via a detailed
study of higher commutators, one can determine the first two stages in a
minimal resolution of the dual of $A$, which also establishes the
isomorphism $A \cong \Omega^2(\bbF_2)$.  The group $V(n)$ corresponds to the
unique non-trivial element in $H^2(G_F^{[2]},\Omega^2(\bbF_2))$;
this element in fact restricts non-trivially on every cyclic
subgroup, whence the extension above is totally non-split.

Now for the \textit{W}-group we have a
extension of elementary abelian groups
\[
1\rightarrow\Phi({\mathcal G}_F) \rightarrow {\mathcal G}_F \rightarrow
G_{F}^{[2]}\rightarrow 1,
\]
where $\Phi$ denotes the Frattini Subgroup of ${\mathcal G}_F$. 
Moreover, elaborating on what we have mentioned this can be expressed as
an extension
\[
1\to U\to G_F^{\{3\}}\to {\mathcal G}_F\to 1
\]
where $U=\Gal(F^{\{3\}}/F^{(3)})$.

\begin{definition}
\label{c-field}
A field whose Witt ring is isomorphic to some group ring of the form
${\mathbb Z}/n\bbZ[G]$ is said to be a \textit{C}-\textsl{field} (here
$n$ is permitted to take the value $0$).
\end{definition}

In \cite[1.9]{Ware} a number of equivalent conditions are given for a
field to be a \textit{C}-field. 
We remark that the $p$-adic field $\bbQ_p$, where $p$ is an
odd prime, is a $C$-field.  Furthermore, $\bbR((t_1)) \cdots
((t_n))$, $\bbC((t_1)) \cdots ((t_n))$, and direct limits of such
fields, are also \textit{C}-fields. Further
examples can be found in \cite[p.~46]{Lam2}.

Now suppose that $K
\subset L$ is a Galois extension of fields, and that there is a
subgroup $N$ of $\dot K / \dot K^2$ such that $L = K[\sqrt{N}]$.  
Again from Kummer theory we know that there is a ``perfect pairing'' 
\[
N \times \Gal(L/K) \rightarrow \mu_2.
\]
Suppose further that $F \subset K$ is another field, and that the
extensions $F \subset K$ and $F \subset L$ are also Galois.  Then
there is a conjugation action of $\Gal(K/F)$ on $\Gal(L/K)$, and
an action of $\Gal(K/F)$ on $N$.  These actions are compatible in the
sense that the following lemma holds, see for example \cite[p.~101]{GMi}
or \cite{Waterhouse}.

\begin{lemma}
If $\sigma \in \Gal(K/F)$, $n \in N$, $\gamma \in \Gal(L/K)$, and $\langle
\cdot , \cdot \rangle$ is the Kummer pairing, then 
\[
\langle \sigma n , \gamma \rangle = \langle n , \gamma^{\sigma^{-1}} \rangle,
\]
i.e. the action of $\Gal(K/F)$ on $\Gal(L/K)$ and $N$ is compatible
with the Kummer pairing. 
\end{lemma}
 
\end{section}

\begin{section}{The \textit{V}-group of a \textit{C}-field}

\label{s:F-is-C-implies-W=V}

In this section we will prove one of our main theorems,
namely

\begin{theorem}
\label{t:C-field-implies}
If $F$ is a field of characteristic different from two,
then the \textit{W}-group of $F$ is the same as the \textit{V}-group of $F$
if and only if $F$ is a \textit{C}-field.
\end{theorem}

We briefly describe the plan of the proof.  To show
that whenever $F$ is a \textit{C}-field we have ${\mathcal G}_F =
G_F^{\{3\}}$, we use a classification of \textit{C}-fields
\cite[1.1]{Ware} and work case-by-case to show that $F^{\{3\}} =
F^{(3)}$. To prove the converse, we note that by \cite[1.9]{Ware}, if
$F$ is not a \textit{C}-field, then there exists a binary 
quadratic form of a certain type over $F$, which implies the
existence of a certain permanent cycle in $E_\infty^{1,1}$, the
mod 2 Lyndon-Hochschild-Serre 
spectral sequence for the group extension $1\to\Phi \rightarrow {\mathcal
G}_F \rightarrow G_F^{[2]}\to 1$.  An identification of $E_\infty^{1,1}$
using Merkurjev's theorem \cite{Mer} then shows that $G_F^{[2]}$ must
act nontrivially on $G_{F^{(2)}}^{[2]}$, and so it follows that
${\mathcal G}_F \neq G_F^{\{3\}}$.  

As mentioned above our proof depends on Ware's classification of
\textit{C}-fields, specifically the following result
(recall that $s(F)$ denotes the level of the field $F$):

\begin{proposition}
Let $F$ be a \textit{C}-field.  If $F$ is formally real, then $F$ is
superpythagorean, while if $F$ is not formally real, $\level(F) = 1
\text{ or } 2$. 
\end{proposition}

\begin{remark}
A field $F$ is \textsl{pythagorean} if $F^2 + F^2 = F^2$.  A formally real field
$F$ is \textsl{superpythagorean} if it is pythagorean and given any subgroup $S$
of index $2$ in $\dot{F}$ such that $-1 \notin F$, $S$ is an ordering of $F$.  For further
details, see \cite[p.~44]{Lam2}.
\end{remark}

It is clear from the proposition above that we can prove
one implication in Theorem~\ref{t:C-field-implies} by studying three cases: the
superpythagorean case, where $\level(F) = \infty$, the case $\level(F)
= 1$, and the case $\level(F) = 2$.  Indeed, it is a consequence
of the following three lemmas. 

\begin{lemma}
\label{l:superpythagorean}
If $F$ is a superpythagorean field then
${\mathcal G}_F = G_F^{\{3\}}$.
\end{lemma}

\begin{lemma}
\label{l:level1}
If $F$ is a \textit{C}-field with $\level(F) = 1$ then ${\mathcal G}_F =
G_F^{\{3\}}$. 
\end{lemma}

\begin{lemma}
\label{l:level2}
If $F$ is a \textit{C}-field with $\level(F) = 2$ then
${\mathcal G}_F = G_F^{\{3\}}$.
\end{lemma}

\begin{remark}
The proofs of Lemmas~\ref{l:superpythagorean}, \ref{l:level1}, and
\ref{l:level2} borrow heavily from the study of the \textit{W}-groups of
the relevant types of fields in \cite{MS2}.  In that paper, the
possible \textit{W}-groups of a \textit{C}-field were determined, and
their structure turns out to depend on the same
sort of structural properties of $\dot F / \dot F^2$ as are used in
the proofs of the following lemmas.   
To clarify things, we provide a list
of the possible \textit{W}-groups of a \textit{C}-field.  
\end{remark}

\begin{proof}[Proof of Lemma~\ref{l:level1}]
First we discuss the possible \textit{W}-groups of \textit{C}-fields
$F$ with $\level(F) = 1$.  If $|\dot{F}/\dot{F}^2| = 1$, $F$ is
quadratically closed, so $F^{\{3\}} = F^{(2)} = F$ and ${\mathcal G}_F
= 1$. If $|\dot{F}/\dot{F}^2| = 2$, then ${\mathcal G}_F =
\mathbb{Z}/4\bbZ$ or $\mathbb{Z}/2\bbZ$; by \cite[3.12.3]{MS2} if
${\mathcal G}_F = \mathbb{Z}/2\bbZ$, then $F$ is euclidean, and in
particular $-1$ is not a square in $F$, so this case is impossible for
a \textit{C}-field of level $1$.  Thus we see that if $F$ is a
\textit{C}-field of level $1$ and $|\dot{F}/\dot{F}^2| = 1 \text{ or }
2$, then ${\mathcal G}_F = 1 \text{ or } \mathbb{Z}/4$, accordingly. 

Now let us suppose that $\dot{F}/\dot{F}^2 = \oplus_I \mathbb{Z}/2\bbZ$,
where $|I| \geq 2$.  Then from \cite[3.13]{MS2} we have that 
the \textit{W}-group of $F$ is $\prod_I \mathbb{Z}/4\bbZ$. 

So we must prove that when ${\mathcal G}_F = \prod_I \mathbb{Z}/4\bbZ$,
$F^{\{3\}} = F^{(3)}$.  As we noted in the first paragraph of this
proof, if $I = \varnothing$, this is trivial, so (although it is not
logically necessary) we can eliminate the need to keep track of
trivial cases by assuming that $|I| \geq 1$.  

To show that $F^{\{3\}} = F^{(3)}$ we will show show for
each $\gamma \in \dot{F}^{(2)}$ that $F^{(2)}(\sqrt{\gamma})/F$ is
Galois, so that $F^{\{3\}} \subset F^{(3)}$, and hence $F^{\{3\}} =
F^{(3)}$.  Our method of demonstrating that 
$F^{(2)}(\sqrt{\gamma})/F$ is Galois will be to prove something a bit
stronger, namely that for each $\gamma \in \dot{F}^{(2)}$ there is an
$a \in \dot{F}$ such that $[\gamma] = [\sqrt{a}]$ in
$\dot{F}^{(2)}/(\dot{F}^{(2)})^2$.  From this it follows that
$F^{(2)}(\sqrt{\gamma}) = F^{(2)}(\sqrt[4]{a})$, which is a Galois
extension of $F$ since $\sqrt{-1} \in F^{(2)}$. 

So, let us take $\gamma \in \dot{F}^{(2)}$ and seek an element 
$a \in \dot{F}$ such that $[\gamma] = [\sqrt{a}]$ in
$\dot{F}^{(2)}/(\dot{F}^{(2)})^2$.  Note first that $\gamma \in
F(\sqrt{a_1}, \ldots, \sqrt{a_n})$ for some elements $a_1, \ldots, a_n$
of $F$, and that we may take $[a_1], \ldots, [a_n]$ to be linearly
independent in $\dot{F}/\dot{F}^2$. Writing $K := F(\sqrt{a_1},
\ldots, \sqrt{a_n})$ and using the proof of \cite[3.13]{MS2}, it
follows that we may complete our partial basis of $\dot{F}/\dot{F}^2$
to a basis ${\mathcal B} = \{[a_1], \ldots, [a_n]\} \cup \{[a_j] \mid j
\in \Omega-\underline{n}\}$ of $\dot{F}/\dot{F}^2$ with the property that 
$\tilde{\mathcal B} = \{[\sqrt{a_1}], \ldots, [\sqrt{a_n}]\} \cup \{[a_j] \mid j
\in \Omega-\underline{n}\}$ is a basis of $\dot{K}/\dot{K}^2$.  (Here
$\Omega$ is an ordinal number, used as an indexing set, and
$\underline{n} := \{1,\ldots,n\}$.) Then
in $\dot{K}/\dot{K}^2$ we have
\[
\gamma = \prod_{i=1}^{n} [\sqrt{a_i}]^{\epsilon_i} \cdot \prod_{j \in
\Omega-\underline{n}} [a_j]^{\epsilon_j},
\]
where $\epsilon_i = 0 \text{ or } 1$ and all but finitely many terms in
the product are $1$, so that  
\[
[\gamma] = \biggl[ \sqrt{\prod_{i=1}^{n} [a_i]^{\epsilon_i} \cdot \prod_{j \in
\Omega-\underline{n}} [a_j]^{2\epsilon_j}} \biggr]
\]
in $\dot{K}/\dot{K}^2$, and hence also in
$\dot{F}^{(2)}/(\dot{F}^{(2)})^2$ since $K \subset F^{(2)}$. 

This completes the proof that $F^{\{3\}} = F^{(3)}$ for
\textit{C}-fields of level $1$.
\end{proof}

The proofs of Lemmas~\ref{l:superpythagorean} and \ref{l:level2} are
quite similar to the proof above and so we give a somewhat
abbreviated presentation.

\begin{proof}[Proof of Lemma~\ref{l:superpythagorean}]
We follow the same plan: we show that for each $\gamma \in
\dot{F}^{(2)}$ there is an $a \in \dot{F}$ such that $[\gamma] =
[\sqrt{a}]$ in $\dot{F}^{(2)}/(\dot{F}^{(2)})^2$.  It follows that 
for any $\gamma \in \dot{F}^{(2)}$ that $F^{(2)}(\sqrt{\gamma})/F$ is
Galois, so that $F^{\{3\}} = F^{(3)}$.  The key point is again that
there is a basis of $\dot{F}/\dot{F}^2$ of the form $\{ [-1] \} \cup
\{[a_i] \mid i \in  I \}$ such that the set $\{[\sqrt{a_i}] \mid i \in
I \}$ is a basis of $\dot{F}^{(2)}/(\dot{F}^{(2)})^2$. 
\begin{remark}
The failure to list $\sqrt{-1}$ as a basis element of
$\dot{F}^{(2)}/(\dot{F}^{(2)})^2$ is not a mistake. If $L$ is a field
of characteristic not equal to 2, then $\sqrt{2} \in L^{(2)}$ and
$\sqrt{-1} \in L^{(2)}$, so that $\sqrt[4]{-1} =
\frac{\sqrt{2}}{2}(1+\sqrt{-1}) \in L^{(2)}$, i.e. $\sqrt{-1}$ is a
square in $L^{(2)}$.
\end{remark}

The proof that $F^{\{3\}} = F^{(3)}$ is then completed my mimicking
the details at the end of the proof of Lemma~\ref{l:level1}.

\end{proof}

The proof of Lemma~\ref{l:level2} will involve an appeal to the
following fact, whose proof can be found in \cite{Ware}.

\begin{lemma}
\label{l:Ware}
Let $K/F$ be a Galois extension and $a$ an element of $\dot{K}$.  Then
the extension $K(\sqrt{a})/F$ is Galois if and only if $\sigma(a) \cdot
a$ is a square in $\dot{K}$ for every $\sigma \in \Gal(K/F)$. 
\end{lemma}

\begin{proof}[Proof of Lemma~\ref{l:level2}]
As in the proof of Lemma~\ref{l:level1}, it is enough to show that for
each $\gamma \in \dot{F}^{(2)}$ that $F^{(2)}(\sqrt{\gamma})/F$ is
Galois.  Since in the case at hand, $F$ is a \textit{C}-field of level
$2$, by \cite[p.~527]{MS2}, we have a basis for
$\dot{F}^{(2)}/(\dot{F}^{(2)})^2$ of the form $\tilde{\mathcal B} =
\{[y],[\sqrt{a_i}] \mid i \in \Omega\}$, where $y \in F(\sqrt{-1})$ and
such that ${\mathcal B} = \{[-1],[a_i] \mid i \in \Omega\}$ is a basis
of $\dot{F}/\dot{F}^2$. 

We claim that in order to show that for
every $\gamma \in \dot{F}^{(2)}$, $F^{(2)}(\sqrt{\gamma})/F$ is
Galois, it suffices to prove this fact for $\gamma$ in $\tilde{\mathcal
B}$. 
\begin{proof}[Proof of claim]
Suppose $\gamma = \gamma_1 \cdots \gamma_s$, where $\gamma_i \in
\tilde{\mathcal B}$. To show that $F^{(2)}(\sqrt{\gamma})/F$ is
Galois, by Lemma~\ref{l:Ware} it is enough to show that
$\sigma(\gamma)\cdot \gamma$ is a square in $F^{(2)}$ for each $\sigma
\in G_F^{[2]}$.  But if this is true for each $\gamma_i \in 
\tilde{\mathcal B}$, the calculation 
\[
\sigma(\gamma)\cdot \gamma = \sigma(\gamma_1 \cdots \gamma_s) \cdot
\gamma_1 \cdots \gamma_s = \sigma(\gamma_1)\cdot
\gamma_1 \cdots \sigma(\gamma_s)\cdot\gamma_s = x_1^2 \cdots x_s^2
\]
shows that it is also true for $\gamma$.
\renewcommand{\qedsymbol}{}
\end{proof}

So we must only show that for each $\gamma \in \tilde{\mathcal B}$
that $F^{(2)}(\sqrt{\gamma})/F$ is Galois.  If $\gamma = \sqrt{a_i}$,
this follows from the fact that $\sqrt{-1} \in F^{(2)}$, while if we
take $\gamma = y$, we may use Lemma~\ref{l:Ware}, and note that
$\sigma(y)\cdot y$ is either $y^2$, which is obviously a square, or
$N_{F(\sqrt{-1})/F}(y)$, which is an element of $F$ and therefore
a square in $F^{(2)}$.
\end{proof}

For the other half of the proof we will use group cohomology
and a result due to Merkurjev.
Our basic goal is the identification of the
$E_{\infty}^{1,1}$-term of the Lyndon-Hochschild-Serre
spectral sequence for the group
extension 
\[ 1\to\Phi \rightarrow {\mathcal G}_F \rightarrow 
G_F^{[2]}\to 1
\]
in
terms of information in the \textit{V}-group of $F$, $G_F^{\{3\}}$. 
To describe this information, we first introduce some notation.

Recall that there is a group extension 
$1\to A \rightarrow G_F^{\{3\}} \rightarrow G_F^{[2]}\to 1$ where the
Pontrjagin dual of $A$ is isomorphic to
$\dot{F^{(2)}}/(\dot{F^{(2)}})^2$.  This 
elementary abelian $2$-group is a $G_F^{[2]}$-module.

Let $J := \dot{F^{(2)}}/(\dot{F^{(2)}})^2$, regarded as a
$G_F^{[2]}$-module.  The socle series of $J$,
\[
0 = J_0 \subset J_1 \subset J_2 \subset \dots \subset J
\]
has finite length $l$ (i.e. $J_l = J$) if $\dot{F}/\dot{F}^2$ is finite.  
This length is an invariant of the field, which we will write as
$l(F)$.  In this socle series, $J_1= J^{G_F^{[2]}}$, $J_2/J_1=
(J/J_1)^{G_F^{[2]}},\dots, J_{i+1}/J_i=(J/J_i)^{G_F^{[2]}}$, etc. We will loosely
refer to this as the length
of the field $F$.

Now we can state our key result 
\begin{proposition}
\label{p:Mer}
In the mod $2$ Lyndon--Hochschild--Serre
spectral sequence for the group
extension 
$1\to \Phi \rightarrow {\mathcal G}_F \rightarrow
G_F^{[2]}\to 1$,
we have that
$E_{\infty}^{1,1} \approx J_2/J_1$.
\end{proposition}

Note that 
the group extension is one for the \textit{W}-group, but that the
group $J_2/J_1$ comes from the \textit{V}-group.  Thus, this fact is
more subtle than it looks.  In fact its proof involves 
the use of Merkurjev's theorem in an essential way.  

Consider the commutative diagram of 
extensions:
\[
\xymatrix{1 \ar[r] & G^{(2)} \ar[r] \ar[d] &  G_q \ar[r] \ar[d] &
G_F^{[2]} \ar@{=}[d] \ar[r] & 1 \\ 
1 \ar[r] & \Phi \ar[r] & {\mathcal G}_F \ar[r] & G_F^{[2]} \ar[r] & 1. \\
}
\]

We shall denote the spectral sequences for the mod 2
cohomology of
${\mathcal G}_F$ and $G_q$ by $E$ and $\bar{E}$ respectively. 
There is an induced map from $E_2^{1,1}$ to $\bar{E}_2^{1,1}$, for which we
will write $\gamma$. 

\begin{lemma}
$\bar{d}_2^{1,1}$ is injective.
\end{lemma}
\begin{proof}
By Merkurjev's theorem, the inflation map from $H^2(G_F^{[2]})$ to
$H^2(G_q)$ is surjective, so $\bar{E}_\infty^{1,1} = 0$, but this is
just the kernel of $\bar{d}_2^{1,1}$. 
\end{proof}

\begin{lemma}
\label{l:kerg=kerd}
$\ker(\gamma) = \ker(d_2^{1,1})$.
\end{lemma}
\begin{proof}
This follows by considering the following commutative diagram and
applying the previous lemma.
\[
\xymatrix{ \bar{E}_2^{1,1} \ar[r]^{\bar{d}_2^{1,1}} & H^3(G_F^{[2]})
\ar@{=}[d] \\  
E_2^{1,1} \ar[r]^{d_2^{1,1}} \ar[u]^{\gamma} & H^3(G_F^{[2]}) \\
}
\]
\end{proof}

\begin{lemma}
\label{l:kerg=J2J1}
$\ker(\gamma) = J_2/J_1$
\end{lemma}
\begin{proof}
We have $\bar{E}_2^{1,1} = H^1(G_F^{[2]}, H^1(G^{(2)}))$ and
$E_2^{1,1} = H^1(G_F^{[2]}, H^1(\Phi))$.  Since $J_1  \approx H^1(\Phi)$
and $J \approx H^1(G^{(2)})$, the result follows from considering the
long exact sequence in cohomology 
associated to the short exact
sequence of coefficients $J_1 \hookrightarrow J \twoheadrightarrow
J/J_1$, since $\gamma$ is just the map $H^1(G_F^{[2]}, J_1)
\rightarrow H^1(G_F^{[2]}, J)$. 
\end{proof}

Now the proof of Proposition~\ref{p:Mer} is merely a matter of
stringing the lemmas together. 

\begin{proof}[Proof of Proposition~\ref{p:Mer}]
Note that $E_\infty^{1,1} = E_3^{1,1} = \ker(d_2^{1,1})$, as there are
no further differentials in this part of the spectral sequence. By
Lemma~\ref{l:kerg=kerd}, $\ker(d_2^{1,1}) = \ker(\gamma)$, which is
$J_2/J_1$ by Lemma~\ref{l:kerg=J2J1}. 
\end{proof}

We can now show that if $F$ is not a \textit{C}-field, then
the \textit{V}-group of $F$ is \emph{not} the \textit{W}-group of $F$.
The proof of this fact is based on an examination of the spectral
sequence for the cohomology of the \textit{W}-group (\ref{p:E11}) and an
identification of part of the $E_\infty$-term of that spectral
sequence with information contained in the \textit{V}-group
(Lemma~\ref{l:123}). 
 
\begin{lemma}
\label{l:123}
The following conditions are equivalent:
\begin{enumerate}
 \item $J = J_1$
 \item $F^{\{3\}} = F^{(3)}$
 \item $J_2 = J_1$.
\end{enumerate}
\end{lemma}
\begin{proof}

\begin{proof}[$1 \Leftrightarrow 2$] 
We need the fact that an extension $F^{(2)}(\sqrt{\gamma})/F$ is Galois
if and only if $[{\gamma}] \in J_1$.  This follows from
Lemma~\ref{l:Ware}.

We have shown that $F^{(3)}$ is obtained
from $F^{(2)}$ by adjoining 
the square roots of elements in $J_1$, while $F^{\{3\}}$ is obtained 
by adjoining the square roots of elements in $J$.  Therefore we have
$F^{\{3\}} = F^{(3)}$ if and only if $J = J_1$. 
\renewcommand{\qedsymbol}{}
\end{proof}

\begin{proof}[$1 \Leftrightarrow 3$]To show that $J_2 = J_1$ implies $J = J_1$ we note that for any element
$j \in J$, the orbit of $j$ under the action of $G_F^{[2]}$ is finite.  This is
because we may choose a representative of $j$, which is 
an element of $\dot{F}^{(2)}$, hence algebraic over $F$, so that any
image of this element under the action of the Galois group satisfies
the same minimal polynomial.  
Therefore, if $\bar{\jmath} \in J/J_1$, the $G_F^{[2]}$ submodule
generated by $\bar{\jmath}$ is finite and has a fixed point by the usual
counting argument.  Thus if $J/J_1 \neq 0$, $J_2 \neq J_1$, which is
the contrapositive of the desired statement.  The converse is trivial.
\renewcommand{\qedsymbol}{}
\end{proof}

\end{proof}

Now we turn to the study of the spectral sequence for the group
extension
\[
1\to \Phi \rightarrow {\mathcal G}_F \rightarrow
G_F^{[2]}\to 1.
\]
We will show that if $F$ is not a \textit{C}-field, then
there are classes in the spectral sequence that survive to
$E_\infty^{1,1}$.  More precisely we will prove: 

\begin{proposition}
\label{p:E11}
Suppose $F$ is \emph{not} a \textit{C}-field.  Then in the spectral
sequence for the cohomology of ${\mathcal G}_F$, $E_\infty^{1,1} \neq 0$.
\end{proposition}

In the proof of this proposition we will need the following lemma
which follows from \cite[theorem~2.7, p.~58]{Lam1} and
\cite[p.~255]{Pierce}.

\begin{lemma}
\label{l:Lam}
Let $L$ be a field and $a$, $b$ elements of $\dot{L}$.  Writing $[a]$
and $[b]$ for the classes of $a$ and $b$ in $\dot{L}/\dot{L}^2$, or
equivalently in $H^2(\Gal(L_q/L);{\mathbb F}_2)$, we have $[a] \cup [b]
= 0$ if and only if $ax^2+by^2 -z^2 = 0$ has a nontrivial solution
$(x,y,z) \in L^3$. 
\end{lemma}

\begin{remark}
More precisely, we need the fact that this lemma remains true if the
Galois group of the quadratic closure $L_q$ is replaced by the
\textit{W}-group ${\mathcal G}_L$.  But this is true by the
construction of the \textit{W}-group and Merkurjev's theorem. See
\cite[3.14]{AKM}.
\end{remark}

\begin{proof}[Proof of Proposition~\ref{p:E11}]
Since $F$ is not a \textit{C}-field, by \cite[1.9]{Ware} there exists
an anisotropic binary form $B$ over $F$ such that the set $D$ of values of
$B$, regarded as a subgroup of $\dot{F}/\dot{F}^2$, has at least three
elements. (As $\dot{F}/\dot{F}^2$ is an elementary abelian $2$-group,
so are any of its subgroups, so ``$D$ has at least three elements''
immediately implies ``$D$ has at least four elements''.
The fact that
$D$ is a subgroup follows from the identity $(1-ax^2)(1-ay^2) =
(1+axy)^2 -a(x+y)^2$.)  We may assume
that $B = x^2 - a y^2$ for some $a \in \dot{F}$ by an appropriate
transformation.  Because $B$ is anisotropic, we know that $[a] \neq [1]$
in $\dot{F}/\dot{F}^2$.

Since $|D| \geq 4$, there exist elements $a_2$, $a_3$ in $D$ which are
linearly independent over ${\mathbb F}_2$.  The
statement ``$a_2$ is a value of $B$ in $\dot{F}/\dot{F}^2$'' means
that there exist $x$, $y$, and $z$ in $F$ such that
$a_2x^2 + a y^2 -z^2 = 0$.  By
Lemma~\ref{l:Lam} we see that $[a][a_2] = 0 \in H^2({\mathcal G}_F)$,
and similarly that  $[a][a_3] = 0$. 

This information allows us to construct directly a nonzero element in
the spectral sequence for $H^*({\mathcal G}_F)$.  Because the
relations $[a][a_2] = 0$ and 
$[a][a_3] = 0$ exist in 
$H^*({\mathcal G}_F)$ but not in $H^*(G_F^{[2]})$, there exist $z_2$,
$z_3$ in $H^1(\Phi)$ such that $d_2^{0,1}(z_2) = [a][a_2]$ and 
$d_2^{0,1}(z_3) = [a][a_3]$.  Then, setting $\lambda = [a_2] \otimes
z_3 + [a_3] \otimes z_2 \in E_2^{1,1} = H^1(G_F^{[2]}) \otimes
H^1(\Phi)$, we can compute 
\[
d_2^{1,1}(\lambda) = 
[a_2] \cdot d_2^{0,1}(z_3) + [a_3] \cdot d_2^{0,1}(z_2) =
[a_2][a][a_3]+ [a_3][a][a_2] = 0.
\]
Thus $\lambda$ survives to $E_3^{1,1}$ and hence to $E_\infty^{1,1}$
by its position in the spectral sequence. 
\end{proof}

The results above immediately imply the converse implication
in the statement of
\ref{t:C-field-implies} and so the proof is complete.
\end{section} 

\begin{section}{Surjectivity in Merkurjev's Theorem and an Analogue of Hilbert's Theorem 90}

We have seen how methods from group cohomology can be used to characterize
certain fields. In this section we make an explicit analysis of the
kernel of the cohomological
differential $d_2$ used previously and interpret the
result in field-theoretic terms. From this we obtain a result which
is reminiscent of Hilbert's Theorem 90.
 
Let $\{[a_i] \mid i\in I\}$ be a basis for $\dot F/\dot F^2$,
which we identify with $H^1(G_F^{[2]})$. Let $\{\sigma_i
\mid i\in I\}$ be a minimal set of generators of $G_F^{[2]}$ such that
${\sigma_i(\sqrt{a_j})}/{\sqrt{a_j}}= (-1)^{\delta_{i,j}}$, where
$i,j\in I$ and $\delta_{i,j}$ is the Kronecker delta function.
One may view the $a_i$ as a dual basis to the $\sigma_i$. 

We are now able to state the result:

\begin{theorem}
\label{p:J90}
Let $\left\{[j_{i_k}] \mid 1\leq k\leq l\right\}$ be a  finite
set of elements of $J_1$. Then
\begin{equation}
\label{e:condition}
    \sum_{k=0}^l[a_{i_k}]d_2^{0,1}\left([j_{i_k}]\right) =0
\end{equation}
if and only if there exists an element $[j]\in J_2$ such that
\begin{enumerate}
\item $[{\sigma_{i_k}(j)}/{j}] =[j_{i_k}]$ for each $k\in\{1,\ldots,l\}$ and 
\item $[\sigma_u(j)] =[j] $ for each $u \notin \{i_1,\ldots,i_l\}$. 
\end{enumerate}
\end{theorem}

It is worth noting that one implication of the theorem (the ``if'' part)
can be proved directly, without the use of Merkurjev's theorem or spectral
sequences.  We indicate how this can be done in an appendix
(\S\ref{s:appendix}).

\begin{proof}[Proof of Theorem~\ref{p:J90}]
Recall that we have the exact sequence 
\[
     0\to {J_2}/{J_1} \rightarrow  H^1(G_F^{[2]},J_1)\rightarrow
     H^1(G_F^{[2]},J) 
\]
arising from the long exact sequence in cohomology associated to the
short exact sequence of coefficients $J_1 \hookrightarrow J
\twoheadrightarrow J/J_1$. Note that since the action of $G_F^{[2]}$
on the coefficient group $J_1$ is trivial, we may identify
$H^1(G_F^{[2]},J_1)$ with $\Hom(G_F^{[2]},J_1)$. 
In other words ${J_2}/ {J_1} \cong \Ker (\Hom(G_F^{[2]},J_1)
\to H^1(G_F^{[2]},J))$. We shall make this isomorphism explicit. 
Let $[j]\in {J_2}/{J_1}$.
Then we can
associate to this $[j]$ the function 
$f_{[j]}:G_F^{[2]}\to J_1$ which is given by the formula:   
\[
    f_{[j]}(\sigma) = \left[\frac{\sigma(j)}{j}\right]\in J_1\,.
\]
It follows from the definition of socle series that
$[{\sigma(j)}/{j}] \in J_1$.
The connecting homomorphism $J_2/J_1 \rightarrow H^1(G_F^{[2]},J_1)$ is given
by $[j] \mapsto
f_{[j]}$.
In other words, we have a one-to-one correspondence $J_2/J_1 \rightarrow 
\Ker(\Hom(G_F^{[2]},J_1) \to H^1(G_F^{[2]},J))$ given by $[j] \mapsto
f_{[j]}$. 
Using the basis $\{a_i\}$ introduced earlier
each function $f_{[j]}$ can be written as 
\[
    f_{[j]} = [a_{i_1}] \otimes [j_{i_1}] + \cdots +
    [a_{i_k}] \otimes [j_{i_k}] \in H^1(G_F^{[2]})\otimes J_1, 
\]
simply because $\Hom(G_F^{[2]}, J_1) = \Hom(G_F^{[2]}, \bbZ/2\bbZ) \otimes
J_1$. 
Observe that our function $f_{[j]}$ has the values:
\[
     f_{[j]}(\sigma_{i_1}) =[j_{i_1}], \ldots, f_{[j]}(\sigma_{i_k}) =[j_{i_k}]
\]
and
\[
   f_{[j]}(\sigma_l) =[1]\,,\quad \text{for each}\quad l\not=
   i_1,\ldots, i_k\,. 
\]
(Recall that each continuous homomorphism $G_F^{[2]}\to J_1$ has only
finitely many values by definition of the Krull topology on
$G_F^{[2]}$.)

Combining the above
with the fact that $J_2/J_1\equiv
E_\infty^{1,1}$ (\ref{p:Mer})
and the definition of $E_\infty^{1,1}=\ker d_2^{1,1}$ completes the proof.  

\end{proof}

\begin{example}
Let $F$ be a field such that $\mcG_F \cong W(2)$.  In this case
$G_q$ is the free pro-2-group on two generators, and 
as we saw before, the \textit{V}-group is a quotient
expressible as an associated ``abelianized extension''.

We can choose a basis $\{[a_1],[a_2]\}$ of ${\dot F}/{\dot F^2}$
and we can choose as our generating $k$-invariants of $W(2)$ the
following elements of 
$H^2(G_F^{[2]})$: 
\[
        q_1= d_2^{0,1}([j_1])=[a_1][a_1], \quad    q_2=
        d_2^{0,1}([j_2])=[a_2][a_2], \quad  
        q_3 = d_2^{0,1}([j_3])=[a_1][a_2].
\]
Thus $J_1=\langle[j_1],[j_2],[j_3]\rangle \subset J$. 

Set also $\lambda_1=[a_1]\otimes [j_3]+[a_2]\otimes[j_1]\in E_2^{1,1}$,
$\lambda_2=[a_1]\otimes[j_2]+[a_2]\otimes[j_3]\in E_2^{1,1}$.

We see that $\langle \lambda_1,\lambda_2\rangle=\ker d_2^{1,1}$ and
therefore $\lambda_1,\lambda_2$ form a basis of $E_\infty^{1,1}$. We
will identify $\lambda_1,\lambda_2$ with elements
$[\lambda_1], [\lambda_2]\in J_2/J_1$ (here we are using
\ref{p:Mer}). 

We pick as usual generators $\sigma_1,\sigma_2$ for $G_F^{[2]}$. 
Then we have $\sigma_1([\lambda_1])=[\lambda_1][j_3]$ and 
$\sigma_2([\lambda_1])=[\lambda_1][j_1]$ similarly the expression for
$\lambda_2$ as an element of $E^{1,1}$ determines
the action of $G_F^{[2]}$ on $[\lambda_2]$:
$\sigma_2[\lambda_2]=[j_3][\lambda_2]$ and
$\sigma_1[\lambda_2]=[\lambda_2][j_2]$.  

It is convenient to picture our $G_F^{[2]}$-module $J_2$ as follows:
\[
\xymatrix{ & [\lambda_1] \ar[rd]^{\Sigma_1} \ar[ld]_{\Sigma_2} &  &
[\lambda_2] \ar[rd]^{\Sigma_1} \ar[ld]_{\Sigma_2} \\
[j_1] & & [j_3] & & [j_2] \\
}
\]
Here $\Sigma_2[\lambda_2]:=
(\sigma_2-1)[\lambda_2]:={\sigma_2[\lambda_2]}/{[\lambda_2]}$, 
etc. 

As mentioned above, $G_q$ is the
free pro-$2$-group on two generators (see \cite[4.3,
pp.~33,34]{S-GC}). Recall that the
well-known formula, due to Schreier, on the number of any minimal set
of generators of any subgroup of finite index of a free group has a
pro-$p$-analogue \cite[p.~49, 6.3]{Koch}.  
Set $G^{(2)}:= \Gal(F_q/F^{(2)})$. Then $G^{(2)}$ is an open subgroup
of $G_q$, since here $|\dot F / \dot F^2| = 4 < \infty$. Then
by the aforementioned formula we have  
\[
   \dim_{\Fbb_2}H^1(G^{(2)}) = 4(\dim H^1(G_q)-1) +1. 
\]
Hence we see that $\dim_{\Fbb_2} H^1(G^{(2)}) = 5$. From Pontrjagin
duality we know that
\[
J={\dot F^{(2)}}/{(\dot F^{(2)})^2}\cong
H^1(G^{(2)}).
\]
Thus $\dim_{\Fbb_2}J =5$.
Since $\dim_{\Fbb_2}J_1 = 3$ and $\dim_{\Fbb_2}J_2/J_1=2$ we see that
$l(J) =2$. This means $J_2 = J$.

Therefore we have the following extension
\[
   1\to A\to G_F^{\{3\}}\to \bbZ/2\bbZ \times \bbZ/2\bbZ\to 1
\]
Our group $A \cong (\bbZ/2\bbZ)^5 = \langle\hat\sigma_1^2,
\hat\sigma_2^2, [\hat\sigma_1, \hat\sigma_2],
[\hat\sigma_1^2,\hat\sigma_2], [\hat\sigma_2^2,\hat\sigma_1]
\rangle$. 
$J$ is just the dual of $A$. This means $J \cong  H^1(A)
\cong \Omega^{-2}(\bbF_2)$. Observe also
that our basis  $\{[j_1],[j_2],[j_3],[\lambda_1],[\lambda_2]\}$ is
dual to the basis of $A$ given above with respect to the
Kummer pairing. Thus we may write  
\[
   [j_1]=(\hat\sigma_1^2)^*, \quad [j_2]=(\hat\sigma_2^2)^*,  \quad[j_3]=([\sigma_1,\sigma_2])^*\,,
\]
\[
     [\lambda_1]= ([\sigma_1^2,\sigma_2])^*,\quad [\lambda_2]=([\sigma_2^2,\sigma_1])^*\,.
\]
It is not hard to use the description above to make explicit
the existence of a  central extension
$1 \ra \bbZ/2\bbZ \times \bbZ/2\bbZ \ra  V(2) \to W(2) \to 1$.
   
\end{example}

\end{section}

\begin{section}{Length and Cohomological Dimension}

We turn now to the case where $\dot F / \dot F ^2$ is finite, so the
order of this group is $2^n$ for some $n$.  In this case $J$ is a
finite-dimensional $G_F^{[2]}$-module.  Therefore its socle series 
$J_1 \subset J_2 \subset \cdots \subset J_l$ has length $l$ for some
$l \in \bbN$.  This number is an interesting field-theoretic
invariant and there are conjectural connections between $l(F)$ and the
cohomological dimension of $G_q$. More generally, for
$G_F^{[2]}$-modules $M$, we will write $l(M)$ for the length of the
socle series of $M$, provided that this length is finite.

If $G$ is a profinite group, its cohomological dimension
$\cd(G)$ is defined as $\sup \cd_p(G)$ where $\cd_p(G)$ is the smallest integer
$n$
such that the $p$-primary component of $H^q(G,A)$ is zero for
all discrete, torsion $G$-modules $A$ and all integers $q>n$
(see \cite[3.1]{S-GC}). Recall that if $G$ is a pro-$2$-group then
$\cd(G) = \cd_2(G) \leq n$ iff $H^{n+1}(G,\bbF_2) = 0$ (see
\cite[p.~27]{S-GC}).

One of the basic connections we wish to consider
appears in 
\cite[6.2.8]{G}: 

\begin{question}
\label{q:Gao}
For which non-formally real fields $F$ with $|\dot F / \dot F ^2|
= 2^n$ is it true that $l(F) + \cd(G_q) = n+1$?
\end{question}

We shall provide some examples of classes of fields for which equality
holds in the equation above.
In the remainder of this section we show that equality holds in the
following cases: first, if $F$ (in addition to satisfying the hypotheses of the
question) is also a  $C$-field, second, if $\mcG_F \approx W(n)$, and
third, if $F$ is a local field. 

\begin{proposition}
\label{p:F=C-field}
Let $F$ be a $C$-field which is not formally real, and suppose in
addition that $|\dot F / \dot F ^2| = 2^n$. Then $l(F) + \cd(G_q) =
n+1$. 
\end{proposition}
\begin{proof}
We have shown that  for $C$-fields, $\mcG_F =
G_F^{\{3\}}$, so that $J = J_1$ and $l(F) =1$.  So we must show that
$\cd(G_q) = n$. From the isomorphism between Galois cohomology and the
Witt ring given by the Milnor conjecture, it is enough to show that
$I^n F/ I^{n+1} F  \neq 0$ and that $I^{n+1} F/ I^{n+2} F  = 0$, where
$I^n F$ is the $n$-th power of the fundamental ideal in the Witt
ring.
Ware \cite{Ware} has shown for fields satisfying our hypotheses that
$I^nF \neq 0$, while it follows from a result of Kneser
\cite[p.~317]{Lam1} and the basic theory of Pfister forms that
$I^{n+1}F  = 0$ in this case.
\end{proof}

We shall show in what follows that the $G_F^{[2]}$-module $J$ can be
studied in a very concrete way.  In the following examples we will use
the identification of $J_2/J_1$ as $E_\infty^{1,1}$ to calculate
$\dim_{\bbF_2}(J_2/J_1)$ in some interesting special cases. 

\begin{example}
\label{e:GF=Wn}
Suppose that $\mcG_F = W(n)$.  Then it follows from  \cite[6.12]{AKM}
that $\dim_{\bbF_2}(E_\infty^{1,1}) = n(n+1)(n-1)/3$.
In fact in this case the group $G_q$ is a free pro-2-group, 
which means that it has cohomological dimension equal to one.
In this case we have identified the module $J$ with the $G_F^{[2]}$
module $\Omega^{-2}(\bbF_2)$, which is known to have a socle series
of length $n$. Hence $l(F) + \cd(G_q)=n + 1$, as asserted above.
\end{example}

\begin{remark}
To see that the $l(\Omega^{-2}2(\bbF_2)) = n$, one can use
Remark~\ref{r:soc-ser} to see that it is enough to prove the
surjectivity of $N_{F^{(2)}/F}$.  This surjectivity follows by
repeatedly applying surjectivity for quadratic extensions $K$ which connect
$F$ to $F^{(2)}$.  The surjectivity for quadratic extensions follows
from the fact that each quaternion algebra defined over $K$ will split
over $K$.  To see that this last assertion is true, note that $G_q(F)$
is free, and hence that $G_q(K) \subset G_q(F)$ is also free.
\end{remark}

It seems natural to explore the situation for 2-dimensional groups.
We recall

\begin{definition}
A Demu\v skin group (at the prime $2$) is a pro-$2$-group $G$ 
which is a two-dimensional Poincar\'e duality group,
i.e. $H^2(G,\bbF_2) = \bbF_2$, and the cup product
\[
H^1(G,\bbF_2) \times H^1(G,\bbF_2) \to H^2(G,\bbF_2)
\]
is a perfect
pairing. 
\end{definition}

\begin{example}
\label{e:J2J1}
Suppose that $F$ is a local field and $|\dot F / \dot F ^2| = 2^n$.
Then $G_q$ is a Demu\v skin group (see \cite[4.5]{S-GC}, 
and in particular we have $\cd G_q = 2$ and $\dim_{\bbF_2} H^2(G_q) =
1$.  It follows from this and \cite[3.12]{AKM} that in the 
extension $1\to\Phi \rightarrow \mcG_F \rightarrow
G_F^{[2]}\to 1$ we have  $\dim_{\bbF_2} \Phi = \binom{n+1}{2}-1$. We can
now calculate $\dim_{\bbF_2} E_2^{1,1} = n \cdot (\binom{n+1}{2}-1)$ and
$\dim_{\bbF_2} E_2^{3,0} = n(n+1)(n+2)/6$; since $H^3(G_q) = 0$, it
follows that $d_2^{1,1}$ is surjective and that
$\dim_{\bbF_2}(E_\infty^{1,1}) = n(n-2)(n+2)/3 = \dim_{\bbF_2} J_2/J_1$.
To make this example more concrete, we note that this calculation
of course applies to the special cases $F = \bbQ_p$ ($p$ an odd
prime), $F = \bbQ_2$, and $F = \bbQ_2(\sqrt 2)$. 
In the case $F = \bbQ_p$, our calculations give $\dim_{\bbF_2} J_2/J_1
= 0$, i.e. $J_1 = J_2$.  Another way to phrase this fact is to say
that $\bbQ_p$ is a $C$-field. In the case $F = \bbQ_2$, we obtain
$\dim_{\bbF_2} J_2/J_1 = 5$, while in the case $F = \bbQ_2(\sqrt 2)$,
$\dim_{\bbF_2} J_2/J_1 = 16$. 

\end{example}

We have thus calculated $|J_1|=|\Phi^*|$ (here $\Phi^*$ is the
Pontrjagin dual of $\Phi$) as well as $|J_2/J_1|$ for the case
of a local field; in fact our calculations are valid for any field $F$
such that $G_q$ is a Demu\v skin group. 

Let $G$ be a Demu\v skin group and let $H$
be an open subgroup of $G$. Let $r_G=\dim_{\Fbb_2} H^1(G)$ and $r_H
=\dim_{\Fbb_2} H^1(H)$. 
Then one can verify (see \cite[p.~44, ex.~6]{S-GC}) that
\[
    r_H-2 =[G:H](r_G-2)\,.
\]
Conversely, this property characterises Demu\v skin groups
(see \cite{DL}).
We use this characterization to prove

\begin{proposition} If $G_q$ is a  Demu\v skin group and $|{\dot
F}/{\dot F}^2|=2^n$, then  
\[
\dim_{\Fbb_2} J= 2^n(n-2)+2
\]
\end{proposition}
\begin{proof}
Suppose that $|{\dot F}/{\dot F}^2|=2^n$. Let $H=G_F^{(2)} := \Gal
({F_q}/{F^{(2)}})$. Then $H$ is an open subgroup of $G_q$ and
$[G_q:H]=2^n$. Therefore from the above
we obtain $r_H=2^n(n-2)+2$.  
However from Kummer theory we know that 
$r_H=\dim_{\Fbb_2}\left({\dot F^{(2)}}/{(\dot F^{(2)})^2}
\right)=\dim_{\Fbb_2} J$.  
\end{proof}

To further our understanding of the $G_F^{[2]}$-module $J$ we look
again at some of the  examples above:

\begin{example}
Let $F={\mathbb Q}_2$. Then ${\dot F}/{\dot F}^2= \langle [-1], [2],
[5] \rangle$, so  $n=3$ and $ \dim_{\Fbb_2} J = 2^3(3-2)+2 = 10$. 
Notice that this is the same as $\dim_{\Fbb_2}\Phi +\dim_{\Fbb_2}
J_2/J_1 = 6-1 +5$.  In particular $l(\bbQ_2) = 2$, and so the desired
equality holds in this case. 
\end{example}

\begin{example}
Let  $F={\mathbb Q}_2(\sqrt 2)$, so  $n=4$. Then $\dim_{\Fbb_2} J
=16\cdot 2+2 =34$. 
We know already \ref{e:J2J1} that $\dim_{\Fbb_2} \Phi =\binom{5}{2}-1 =9$
and that $\dim_{\Fbb_2} J_2/J_1 = 16$.  If the equality holds in this case,
then we have $l(F) = 5-2 =3$, so that  
$J = J_3$ and $\dim_{\Fbb_2}{J_3}/{J_2} =\dim_{\Fbb_2}{J}/{J_2}
=34-25=9$. 
\end{example}

\begin{example}
Suppose that $F$ is a local field with $|\dot F/\dot F^2| = 2^n$. Then
we can compute $\dim_{\Fbb_2}({J}/{J_2} )$:
\begin{eqnarray*}
       \dim_{\Fbb_2}({J}/{J_2} )
   &=& \dim_{\Fbb_2} J -\dim_{\Fbb_2} J_2/J_1 -\dim_{\Fbb_2} J_1\\
   &=& 2^n(n-2) +2 -\frac{n(n-2)(n+2)}{3} - \binom{n+1}{2} + 1 \\
   &=& 2^n(n-2) + 3-\frac{n(n-2)(n+2)}{3}-\binom{n+1}{2}.
\end{eqnarray*}
\end{example}
We will now make a more in-depth analysis of the module $J$. 
We recall the following
basic fact from 
local class field theory:

\begin{lemma}
Let $F$ be a local field such that $|\dot F/\dot F^2|=2^n$. If 
$N_{F^{(2)}/F}:\dot F^{(2)}\to \dot F$ is the usual norm map then its image
lies in $\dot F^2$. If $F\subset K\subset F^{(2)}$ is a proper quadratic
extension of $F$, then the image of $N_{F^{(2)}/K}:\dot F^{(2)}\to \dot K$
does not lie in $\dot K^2$.
\end{lemma}

\begin{proof}
 From local class field theory (see \cite{CF}) we know that
there is a natural isomorphism
$\Gal(F^{(2)}/L)\cong \dot L /\Image N_{F^{(2)}/L}$ for any intermediate
extension $F\subset L\subset F^{(2)}$. If $L=F$ we obtain the
first statement. For $L=K$ we observe 
(see \cite{Lam1}, page 202) that for any proper quadratic
extension $K$ of $F$, $|\dot K/\dot K^2|\ge |\dot F/\dot F^2|$; however
$|\Gal(F^{(2)}/K)|<|\Gal(F^{(2)}/F)|=|G_F^{[2]}|$. We conclude that
$\Image N_{F^{(2)}/K}$ cannot be contained in $\dot K^2$.
\end{proof}

This relates to our length invariant via the following 

\begin{proposition}
\label{p:J_n-1}
Let $F$ be a field such that $|\dot F/\dot F^2|=2^n$.
Then, in the socle series $J_1\subset J_2\subset\dots\subset J_{l(F)}$,
the submodule $J_{n-1}$ is equal to the kernel of the homomorphism
$N: J\to \dot F/\dot F^2$ induced by the norm map above.
\end{proposition}

\begin{remark}
\label{r:soc-ser}
The proposition above can be extended to a characterization of the
complete socle series $\{J_i\}$; a proof by induction, using
\cite[p.~133]{ELW}, which we omit,
shows that $J_i = \{ [j] \in J \mid N_{F^{(2)}/L}(j) = [1] \in
\dot{L}/\dot{L}^2 \text{ for all } L \text{ with } [F^{(2)}:L] = 2^{i+1}
\}$.
\end{remark}
\begin{proof}[Proof of Proposition~\ref{p:J_n-1}]
Let $F\subset L\subset F^{(2)}$ denote an extension of $F$ such that
$[F^{(2)}:F]=2^k$. Consider $\Gal(F^{(2)}/L)$; then by restriction
$J$ will also be a module over this group. Denote the usual module-theoretic
norm map by $T_L:J\to J$.
 From the definition of the socle series for $J$ it follows
that 
$J_{k-1} = \cap \ker T_L$,
where the intersection is taken over all extensions as above, of co-degree
$2^k$. Let $L$ denote a proper quadratic extension of $F$. The norm map
$T_L$ is induced from the composition of the field-theoretic norm map
$N_{F^{(2)}/L}\colon \dot F^{(2)}\to \dot L$ with the maps
$c_L\colon\dot L\to \dot L/\dot L^2$ and $i_{L/F}\colon\dot L/\dot L^2\to
J$ (the latter induced 
by the inclusion $L\subset F^{(2)}$). Hence if $[j]\in J$, we see
that $[j]\in \ker T_L$ if and only if $c_L(N_{F^{(2)}/L}(j))\in \ker i_{L/F}$.
However from Kummer theory (see \cite{AT}, page 21, theorem 3) we have that
$\ker i_{L/F}=\dot F\dot L^2/ \dot L^2$, and so $T_L(j)=0$ if and only if
$N_{F^{(2)}/L} (j)\in \dot F\dot L^2$. However, from \cite{Lam1}, theorem 3.4, pp. 202,
we know that $\mu\in \dot L$ belongs to $\dot F\dot L^2$ if and only if
$N_{L/F}(\mu)\in \dot F^2$. Therefore, from the transitivity property
of the norm we see that $T_L([j])=0$ if and only if 
$N_{F^{(2)}/F}(j)\in\dot F^2$.
\end{proof}

 From this (\ref{p:J_n-1}) we obtain
\begin{corollary}
Let $F$ denote a local field with $|\dot F/\dot F^2|=2^n$; then
$l(F)=n-1$.
\end{corollary}
\begin{proof}
 From our description of $J_{n-1}$ and the triviality of the map
induced by the norm, we
infer that $J=J_{n-1}$, whence $l(F)\le n-1$. Taking a quadratic
extension $F\subset K\subset F^{(2)}$ we see that in the socle series
for $J$ as a $\Gal(F^{(2)}/K)$-module, $J_{n-2}\ne J$. Hence $l(F) = n
-1$.
\end{proof}

This can be restated as 
\begin{theorem}
\label{t:F=local}
Let $F$ denote a local field with $|\dot F/\dot F^2|=2^n$. Then
$l(F)+ \cd(G_q) =n+1.$
\end{theorem}

It seems rather complicated to verify the relationship above for
other types of fields. We briefly outline a more cohomological
approach.

Given any finitely generated
$\bbF_2 G_F^{[2]}$-module $M$, there is a minimal power $I^t$ of the
augmentation ideal (which is nilpotent) such that $I^tM=0$; $t=\lambda (M)$ 
is called the
Loewy length of $M$. Note that if $M\ne 0$, then $1\le \lambda (M)\le n+1$.
Moreover $M$ is a trivial $G_F^{[2]}$-module if and only if $\lambda (M)=1$.
The length of the socle series of a module is equal to the Loewy length
of its dual.
Consider the following
result, which follows from a theorem due to G. Carlsson (\cite{Ca}):

\begin{proposition}
Suppose that $H$ and $G$ are pro-$2$-groups, which are topologically
finitely generated,  $G$ of finite (continuous) cohomological dimension $k$ at the
prime $2$, and which have finite total mod $2$ cohomology.
Assume in addition that $E$ is an elementary
abelian $2$-group of rank equal to $n$, and that
\[
1 \rightarrow H \rightarrow G \rightarrow E \rightarrow 1
\]
is an extension.  Then it follows that
\[
\lambda (H^1(H,\bbF_2)+\dots +\lambda (H^k(H,\bbF_2))\ge n.
\]
In particular if $H$ and $G$ are $2$-dimensional Poincar\'e duality groups,
$H^1(H,\bbF_2)$ is self-dual,
$\lambda (H^1(H,\bbF_2)) \ge n-1$, and therefore $l(H^1(H, \bbF_2)) \geq
n-1$.
\end{proposition}

As an immediate consequence of the above we obtain a different proof
of the inequality $l(F)\ge n-1$ for fields
$F$ such that $G_q$ is Demu\v skin.
As indicated by this cohomological method, 
the filtration lengths of higher cohomology groups will probably
play a role in any generalization of this question
to groups of larger cohomological
dimension.
\end{section}

\begin{section}{The Cohomology of $V(2)$}

In this section we study the cohomology of the universal
\textit{V}-group $V(2)$ defined previously.  
We will first describe some basic facts about $V(n)$ expressing it in terms
of certain group extensions which will be useful in computing and
interpreting cohomology.
Our point of view will closely
follow the analysis of \textit{W}-groups made in \cite{AKM}, so we
will be brief.  Many of the calculations described below can be done
partially or entirely using a computer algebra system such as
\textsc{Magma} \cite{Magma}, so almost all detailed justifications are omitted.

We begin by taking the standard surjection from 
the free group on $n$-generators $F_n$ onto the elementary abelian
group $E_n=(\bbZ/2\bbZ)^n$, with kernel the free
group on $2^n(n-1)+1$ generators. Associated to this we have the
free abelianized extension
$1\to M\to X(n)\to E_n\to 1$
where $M$ is a $\bbZ E_n$ lattice of rank as above. It is elementary
to verify (see \cite{Ku} for details) that, as a $\bbZ E_n$--module,
$M$ is isomorphic to $\Omega^2(\bbZ)$,
the second dimension shift of the trivial module. Factoring out the submodule
$2M$, we recover an expression for $V(n)$ as an extension of finite groups,
$1\to M/2M\to V(n)\to E_n\to 1$. This has a simple interpretation: the
module $M/2M\cong \Omega^2(\bbF_2)$ has a unique non zero class in its
second cohomology group $H^2(E_n,\Omega^2(\bbF_2))=\bbF_2$; the group
$V(n)$ realizes this extension class. Similarly, the group $X(n)$ is
a Bieberbach group (see \cite{Wolf}) corresponding
to the canonical generator in $H^2(E_n,\Omega^2(\bbZ))\cong\bbZ/|E_n|$; it
is also expressed as an extension
$1\to 2M\to X(n)\to V(n)\to 1$. 

The comments above show that the group $V(n)$
acts freely on a $2^n(n-1)+1$-dimensional torus, with quotient the
classifying space of $X(n)$. In turn this space can be obtained as
the orbit space of a free $E_n$-action on such a torus. The cohomology of these
euclidean space forms is not easy to compute, and we shall see that even the
case $n=2$ poses some interesting technical problems. 

We should also point
out that from the definition of the \textit{V}-groups we can also express it as
an extension
$1\to U_r\to V(n)\to W(n)\to 1$ where $U_r$ is a subgroup
isomorphic to $(\bbZ/2\bbZ)^r$, $r=2^n(n-1)-n-\binom{n}{2}+1$. An important
thing to note is that the rank of the largest elementary abelian subgroup
in $V(n)$ is precisely equal $r + n + \binom{n}{2}$, or in other words
the rank of $V(n)$ differs from the rank of $W(n)$ exactly by
the quantity $r$.

We will now concentrate on the case $n=2$. Our first result is about
the group $X(2)$:
we study its mod $2$ cohomology via the
spectral sequence for the extension $1\to \bbZ^5 \to X(2) \to \bbZ/2\bbZ
\times \bbZ/2\bbZ\to 1$.  The $\bbZ/2\bbZ \times \bbZ/2\bbZ$-module structure of the
kernel is known, so we can obtain the $\bbZ/2\bbZ \times \bbZ/2\bbZ$-module
structure of the cohomology from the fact that $H^*(\bbZ^5, \bbF_2)$
is an exterior algebra.  Let $k=\bbF_2$ denote the coefficient
field, then this exterior algebra can be written as a sum
of indecomposable $k[ \bbZ/2\bbZ \times \bbZ/2\bbZ]$-modules as follows:

\bigskip
\begin{center}
\begin{tabular}{c | c}
Degree & Decomposition \\ \hline
0 & $k$ \\
1 & $\Omega^{-2} k$ \\
2 & $\Omega^{1} k \oplus \Omega^{1} k \oplus F$ \\
3 & $\Omega^{-1} k \oplus \Omega^{-1} k \oplus F$ \\
4 & $\Omega^{2} k$ \\
5 & $k$ \\
\end{tabular}
\end{center}
\bigskip

 From the table it is straightforward to determine the $E_2$-term of
the spectral sequence and even the $d_2$-differential.  It turns out
that there are no further possible differentials, so the spectral
sequence collapses at $E_3$.  The next few paragraphs sketch a
proof of this fact and record some immediate consequences.

The $E_2$-term is not hard to understand as a $H^*(\bbZ/2\bbZ \times
\bbZ/2\bbZ, \bbF_2)$-module since all the modules appearing in the table
are either free, or dimension shifts of the trivial module.
Furthermore, the $d_2^{*,1}$-differential realizes the cohomology
isomorphism $H^*(\bbZ/2\bbZ \times \bbZ/2\bbZ, \Omega^{-2} k) \approx
H^{*+2}(\bbZ/2\bbZ \times \bbZ/2\bbZ,  k)$, since $d_2^{0,1}$ is an
isomorphism by the definition of $X(2)$ (specifically, the extension
class). By duality it
follows that the $d_2^{*,5}$-differential realizes the cohomology
isomorphism $H^*(\bbZ/2\bbZ \times \bbZ/2\bbZ, k) \approx H^{*+2}(\bbZ/2\bbZ
\times \bbZ/2\bbZ, \Omega^2 k)$.  Since we know that the classifying space
for $X(2)$, and hence its cohomology, is finite-dimensional, it
follows that the  $d_2^{*,3}$-differential realizes the cohomology
isomorphism $H^*(\bbZ/2\bbZ \times \bbZ/2\bbZ, \Omega^{-1} k \oplus \Omega^{-1} k)
\approx H^{*+2}(\bbZ/2\bbZ \times \bbZ/2\bbZ, \Omega^{1} k \oplus
\Omega^{1} k)$. 

We now note that all this implies 
\begin{proposition}
In the spectral sequence for the extension $1\to\bbZ^5 \to X(2) \to \bbZ/2\bbZ
\times \bbZ/2\bbZ\to 1 $, $E_3^{p,q} = 0$ for $p>1$, and $E_3 = E_\infty$.
\end{proposition}

A little more attention to detail gives also
\begin{proposition}
The Poincar\'e Series for the cohomology of $X(2)$ is $1 + 2t + 5t^2 +
5t^3 + 2t^4 + t^5$. 
\end{proposition}

A method similar to that of the previous section could be used
to study the cohomology of $V(2)$. However, this group
has order only $2^7$, and its cohomology can thus be studied in great
detail with a computer.  In particular, information on this group is
available at Carlson's well-known web site \cite{C-web}.

As the details of the calculation are complicated, we present only an
outline of the work necessary to determine by hand calculation the
cohomology of $V(2)$.  The main reason for outlining this work is to
point out the existence of a phenomenon which the careful reader will
already have noted in the previous section on $X(2)$.

Now let us sketch the calculation.  Studying the cohomology of $V(2)$
via the spectral sequence for the extension $1\to (\bbZ/2\bbZ)^5 \to V(2)
\to \bbZ/2\bbZ \times \bbZ/2\bbZ\to 1$, we must write the symmetric
algebra of $\Omega^{-2} k$ as a direct sum of indecomposable modules.
If we write $P$ for the direct sum of the three nontrivial permutation
$\bbZ/2\bbZ \times \bbZ/2\bbZ$-modules which are not free, and $F$ for
the free $\bbZ/2\bbZ \times \bbZ/2\bbZ$-module of rank $1$, the
multiplicities of the various modules in the symmetric algebra can be
given in terms of Poincar\'e series: 

\bigskip
\begin{center}
\begin{tabular}{c | c }
Module & Poincar\'e Series \\ \hline
$k$ & $(1-t^4)^{-2}(1+t^3)$ \\
$\Omega^{-2} k$ & $(1-t^4)^{-2}t$\\
$\Omega^{2} k $ & $(1-t^4)^{-2}t^2$\\
$P$ & $(1-t^4)^{-2}(1-t^2)(t^2+t^3+t^4+t^5)$ \\
$F$ & $(1-t^4)^{-2}(1-t)^{-2}[(1-t)^{-1}t^2 + 4t^3 + 4t^5 +
(1-t)^{-1}3t^6]$  \\ 
\end{tabular}
\end{center}
\bigskip

It turns out that, up to projective summands, $S^{2j+1}
\Omega^{-2}k \approx \Omega^{-2} S^{2j} \Omega^{-2}k$, and that the
differential $d_2^{i,2j+1}$ is the natural isomorphism $H^i(\bbZ/2\bbZ
\times \bbZ/2\bbZ, S^{2j+1}) \approx H^{i+2}(\bbZ/2\bbZ
\times \bbZ/2\bbZ, S^{2j})$.  (To study this differential in an ad hoc
manner one can use the restrictions to index $2$ subgroups of $V(2)$---all of
these are isomorphic to $((\bbZ/2\bbZ)^2 \times (\bbZ/2\bbZ)^2)
\rtimes \bbZ/4\bbZ$, 
and these groups are well understood since they are not far from being
wreath products; in fact the cohomology of these groups is detected on
abelian subgroups.) 
As in the previous section, the spectral sequence associated to the
extension collapses at $E_3$ and $E_3^{p,q} = 0$ if $p>1$.

\begin{remark}
\label{r:coincidence}
The interest of the method above lies in the odd coincidences that arise
in the computations---that $\Omega^{-2} S^{2j}$ should be isomorphic
to $S^{2j+1}$, and that the differential $d_2$ should contain the
associated cohomology isomorphism.  
A similar phenomenon for the group $X(2)$ was observed in the
previous section.  It is natural to wonder if these facts can be
obtained by some simple method that does not involve ad hoc
computations.
\end{remark}

In any case, whether from a computer analysis, or following the outline
above, one can obtain a number of facts about the cohomology of
$V(2)$, which are catalogued below.

\begin{proposition}
\label{p:V2}
The Poincar\'e series of $H^*(V(2))$ is 
\[
\frac{1-x+x^2}{(1-x)^3(1-x^2)^2},
\]
which expands to $1 + 2 x + 6 x^2  + 11 x^3  + 22 x^4  + 36 x^5 + 60
x^6 + 90 x^7 + 135 x^8 + 190 x^9 + 266x^{10} +
\cdots$.  
\end{proposition}

\begin{proposition}
$H^*(V(2))$ is detected on abelian subgroups. 
\end{proposition}

\begin{remark}
Note that the proposition above does \emph{not} claim detection on
elementary abelian subgroups.
\end{remark}

We end our analysis of the groups $V(n)$ on a topological note. One of the key
properties of finite $2$-groups with central involutions is that they
act freely on a product of $s$ spheres, where $s$ is the rank of the
largest elementary abelian subgroup. This number is known to be minimal
for this property and interesting results can be derived from this
(see \cite{AY}). One can prove that the universal \textit{V}-groups also satisfy
this, whence they seem to be rather special. In fact one can use the
extension $1\to U_r\to V(n)\to W(n)\to 1$ to prove:

\begin{proposition}
Let $s = dim~\Omega^2(\bbF_2)$ denote the rank of the largest elementary abelian subgroup
in $V(n)$. Then there exist $2^n$-dimensional
real representations $Z_1,\dots ,Z_s$  of $V(n)$ such that its diagonal
action on the product of associated linear spheres
$X=\bbS (Z_1)\times \dots\times\bbS (Z_s)$ is free.
\end{proposition}
\begin{proof}
Let $x_1,\dots ,x_r$ denote a basis for the subgroup $U_r$.
For each $i$ one can construct a 1--dimensional real representation of
the Frattini subgroup of $V(n)$ which restricts non-trivially to
$<x_i>$. This can then be induced to yield a $2^n$-dimensional
representation $Z_i$ for $V(n)$. Doing the analogous construction
for $W(n)$ and pulling back to $V(n)$, we obtain the representaions
$Z_{r+1}, \dots, Z_s$. It is elementary to verify that every involution
in $V(n)$ must acts freely on the product of associated spheres; indeed
the unique maximal elementary abelian subgroup satisfies this by construction.
\end{proof}
\end{section}

\begin{section}{Appendix : An Elementary Proof of Half of
Theorem~\ref{p:J90}}
\label{s:appendix}

In this section we will provide an elementary proof of ``half'' of
Theorem~\ref{p:J90} (i.e. the proposition below), and in doing so
develop a connection between $J_2/J_1$ and certain triple commutators.

\begin{proposition}
\label{p:halfJ90}
If $\left\{[j_{i_k}] \mid 1\leq k\leq l\right\}$ is a  finite
set of elements in $J_1$, and there is an element $[j]\in J$ such that
\begin{enumerate}
\item $[{\sigma_{i_k}(j)}/{j}] =[j_{i_k}]$ for each $k\in\{1,\ldots,l\}$ and 
\item $[\sigma_u(j)] =[j] $ for each $u \notin \{i_1,\ldots,i_l\}$, 
\end{enumerate}
then 
\[
    \sum_{k=0}^l[a_{i_k}]d_2^{0,1}\left([j_{i_k}]\right) =0 \,.
\]
\end{proposition}

\begin{remark}
The observant reader will note that $d_2^{0,1}$ is merely a notation
for a homomorphism determined by our \textit{W}-group, and that we make
no real use of spectral sequences.
\end{remark}

Since this is an alternate proof, for the sake of simplicity we will
assume $|{\dot F}/{\dot F^2}|= 2^n$; this means that we can simplify the
notation by assuming that $\{i_1, \ldots, i_l\} = \{1, \ldots, n\}$.
We leave to the reader the minor
modifications necessary for the case in which ${\dot F}/{\dot F^2}$ is
infinite. Recall from the notational conventions in the
proof of Theorem~\ref{p:J90} that
to each $[j]\in J_2/J_1$
we can associate an
element $\lambda_j\in H^1(G_F^{[2]})\otimes J_1$,
which we shall assume to have the form: $\lambda_j =\sum_{k=1}^n
[a_k]\otimes [j_k]$.
With the notation of this paragraph, equations 1 and 2 of the
proposition above become
$[{\sigma_k(j)}/{j}] = [j_k] \in J_1$
for each $k=1,2,\ldots,n$.
\begin{notation}
$F_j:=F^{(3)}(\sqrt j)$, which is a Galois extension of $F$,
and $G_j := \Gal(F_j/F)$.
\end{notation} 

We will identify the elements of $J_1$ with ``$k$-invariants of
$\mcG_F$'' via $d_2^{0,1}$.

\begin{notation}
\label{n:enters}
In the situation above, if in the expression of $[j]\in
J_1$ as a sum of monomials, the coefficient of $[a_i][a_j]$ is
nonzero, we will say: ``$[a_i][a_j]$ enters the expression of $[j]\in 
J_1$'', or more concisely ``$[a_i][a_j]$ enters $[j]$.''  
\end{notation}

We will now use the Kummer pairing $J_1\times \Phi\to \mu_2$, given by
$\langle[j],\gamma\rangle ={\gamma(\sqrt j)}/{\sqrt j}$.
 From the
discussion of this pairing in \cite[pp.~42--48]{MS1} we 
have: 
\begin{lemma} Let  $1 \leq i \leq k \leq n$, then
\label{l:MS1}
\begin{enumerate}
\item $\langle[j],[\sigma_i,\sigma_k]\rangle = -1$ iff
$[a_i][a_k] \text{ enters } [j]$.  
\item $\langle[j],\sigma_i^2\rangle = -1$ iff $[a_i][a_i] \text{
enters } [j].$ 
\end{enumerate}
\end{lemma}

\begin{notation}
We denote by $\theta$ the nontrivial element of the kernel of quotient map $G_j
\twoheadrightarrow \mcG_F$.
\end{notation}

Let $\hat{\sigma}_k$, $k = 1, \ldots, n$ be extensions of the generators
$\sigma_k$ of $\mcG_F$ to $G_j$. Let $\{\hat\sigma_i,\hat\sigma_k,\hat\sigma_l\}$, $1\leq i<k<l\leq n$,
be a triple consisting of generators of $G_j$.
A direct calculation proves:
                            
\begin{lemma}
\label{l:fact123}
Using the terminology of \ref{n:enters}, we have the following
identities: 
\begin{enumerate}
\item $[\hat\sigma_i, [\hat\sigma_k,\hat\sigma_l]]=\theta$ iff $[a_k][a_l]$
enters $[j_i]$.   
\item $[\hat\sigma_i, \hat\sigma_k^2]=\theta$ iff $[a_k][a_k]$ enters $[j_i]$.
\item $[a_i][a_i]$ does not enter $[j_i]$.
\end{enumerate}
\end{lemma}

Another combinatorial lemma which we will need is:

\begin{lemma} 
\label{l:identities2}
For $\sigma, \tau, \gamma$ elements of $G_j$, we have
the following additional identities:
\begin{enumerate}
\item[1.]
$[\sigma,[\tau,\gamma]][\tau,[\gamma,\sigma]][\gamma,[\sigma,\tau]]=1$.
\item[2.] $[\sigma,\tau]^2=1$ and $[\sigma,\tau]=[\tau,\sigma]$
\item[3.] $[\sigma^2,\tau]=[\sigma,[\sigma,\tau]]=[\tau,\sigma^2]$
\end{enumerate}
\end{lemma}

Observe that part 1. of this lemma is the Jacobi identity, which is
valid in any metabelian group (see \cite{BL}).

\begin{proof}[Proof of Proposition~\ref{p:halfJ90}]
\label{pf:p:halfJ90}
Suppose that $[j]\in J$ is as in the statement of our Proposition.
In order to show that 
$\sum_{\beta =1}^n[a_\beta]d_2^{0,1} ([j_\beta])=0$
it is enough to observe that for each monomial $[a_i][a_k][a_l]$ that
can occur in the expression of the left hand side as an element of
$H^3(G_F^{[2]})$, the corresponding coefficient is zero.
This can be proved using a case by case analysis, in which division into
cases
depends on the multiplicities of the $[a_i]$ which appear in the monomial.
We limit ourselves to providing a complete analysis for the case where
all the multiplicities are $1$; the other cases are very similar and
left to the reader.

Assume that [$1\leq i<k<l\leq n$] then the term $[a_i][a_k][a_l]\in
H^3(G_F^{[2]})$ can occur as the summand of the following terms of our
sum $\sum_{\beta=1}^n[a_\beta]d_2^{0,1}([j_\beta])$: 
\[
    [a_i]d_1^{0,1}[j_i]\,,\qquad [a_k]d_2^{0,1}[j_k]\,\qquad
    \text{and } [a_l]d_2^{0,1}[j_l].
\]    
 From Lemma~\ref{l:fact123}.1 we see that when this term does
occur, one of our triple commutators
$[\hat\sigma_i,[\hat\sigma_k,\hat\sigma_l]$,
$[\hat\sigma_k,[\hat\sigma_i,\hat\sigma_l]]$, or
$[\hat\sigma_l,[\hat\sigma_i,\hat\sigma_k]]$
will be $\theta$ and not the identity.
However from the identity \ref{l:identities2}.1:
\[
  [\hat\sigma_k,[\hat\sigma_i,\hat\sigma_l]]
     [\hat\sigma_i,[\hat\sigma_k,\hat\sigma_l]]
     [\hat\sigma_l,[\hat\sigma_i,\hat\sigma_k]]=1
\]
and the relation $\theta^2=1$ we see that the term
$[a_i][a_k][a_l]$ occurs in our sum
$\sum_{\beta=1}^n[a_\beta]d_2^{0,1}[j_\beta]$ 
either zero or two times.

This completes the analysis of the case when
all the multiplicities are $1$; an analysis of the other cases proves
that all terms $[a_i][a_k][a_l]$ occur an even number of times in the
sum and thus we have
\[
  \sum_{\beta=1}^n [a_\beta]d_2^{0,1}([j_\beta]) =0
\]
as desired.
\end{proof}

\end{section}

\end{document}